\documentclass[12pt]{article}
\usepackage{amssymb,latexsym}
\def\QQ{{\mathbb Q}}
\def\RR{{\mathbb R}}

\newtheorem{formula}{}[section]

\newtheorem{definition}[formula]{\indent Definition}
\newtheorem{corollary}[formula]{\indent Corollary}
\newtheorem{remark}[formula]{\indent Remark}
\newtheorem{lemma}[formula]{\indent Lemma}
\newtheorem{theorem}[formula]{\indent Theorem}

\def\thrm{\begin{theorem}}
\def\thrml#1{\begin{theorem}\label{#1}}
\def\ethrm{\end{theorem}}
\def\rmrk{\begin{remark}}
\def\rmrkl#1{\begin{remark}\label{#1}}
\def\ermrk{\end{remark}}
\def\dfntn{\begin{definition}}
\def\dfntnl#1{\begin{definition}\label{#1}}
\def\edfntn{\end{definition}}
\def\nmrt{\begin{enumerate}}
\def\enmrt{\end{enumerate}}

\def\qtn{\begin{equation}}
\def\qtnl#1{\begin{equation}\label{#1}}
\def\eqtn{\end{equation}}
\def\lmm{\begin{lemma}}
\def\lmml#1{\begin{lemma}\label{#1}}
\def\elmm{\end{lemma}}
\def\crllr{\begin{corollary}}
\def\crllrl#1{\begin{corollary}\label{#1}}
\def\ecrllr{\end{corollary}}

\begin{document}
\title{}
\date{}
\maketitle
\vspace{-0,1cm} \centerline{\bf Tropical Newton-Puiseux polynomials II}
\vspace{7mm}
\author{
\centerline{Dima Grigoriev}
\vspace{3mm}
\centerline{CNRS, Math\'ematique, Universit\'e de Lille, Villeneuve
d'Ascq, 59655, France} \vspace{1mm} \centerline{e-mail:\
dmitry.grigoryev@univ-lille.fr } \vspace{1mm}
\centerline{URL:\ http://en.wikipedia.org/wiki/Dima\_Grigoriev} }

\begin{abstract}
Tropical Newton-Puiseux polynomials defined as piece-wise linear functions with rational coefficients at the variables, play a role of tropical algebraic functions. We provide explicit formulas for tropical Newton-Puiseux polynomials being the tropical zeroes of a univariate tropical polynomial with parametric coefficients.
\end{abstract}

{\bf keywords}: tropical Newton-Puiseux polynomial, zeroes of tropical parametric polynomials


\section*{Introduction}
One can find basic concepts of tropical algebra in \cite{MS}.

Given a tropical univariate polynomial
\begin{eqnarray}\label{1}
f=\min_{0\le k\le n} \{x_k+kY\}
\end{eqnarray}
its tropical zero $y\in \RR$ is such that the minimum in (\ref{1}) is attained for at least two different values of $0\le k\le n$. In this paper we treat the coefficients $X=(x_0,\dots,x_n)$ as parameters and find the zeroes of $f$ as a function in $x_0,\dots,x_n$.

We show that $f$ has exactly $n$ such parametric zeroes $g_1,\dots,g_n$. Each $g_k,\, 1\le k\le n$ is a {\it tropical Newton-Puiseux polynomial}, i.~e. a piece-wise linear function with rational coefficients at the variables. One can represent any tropical Newton-Puiseux polynomial in the form
$$\min_I \{a_I+(I,\, X)\}-\min_J\{b_J+(J,\, X)\}$$
\noindent of a difference (so, the tropical quotient) of two concave piece-wise linear functions, where $I,\, J\in \QQ^{n+1};\, a_I,\, b_J\in \RR$, and $(I,\, X)$ denotes the inner product (cf. \cite{O}).

Tropical Newton-Puiseux polynomials play a role of algebraic functions in the tropical setting. Similar to Newton-Puiseux series in the classical algebra, $g_k(x_0,\dots,x_n)$ provides a tropical zero of $f$ for any point $(x_0,\dots,x_n)\in \RR^{n+1}$. We note that in the classical algebra one considers Newton-Puiseux series just in a single variable, while an advantage of the tropical algebra is that one considers tropical Newton-Puiseux polynomials in several variables.

Observe that if one considers a univariate tropical polynomial with its coefficients being tropical Newton-Puiseux polynomials then in its turn, the tropical zeroes of this tropical polynomial are again tropical Newton-Puiseux polynomials. Thus, one can view the semi-field of tropical Newton-Puiseux polynomials as a tropical algebraic closure of the semi-ring of tropical polynomials.

\section{Tropical Newton-Puiseux polynomials as tropical zeroes}\label{one}

We say that a tropical Newton-Puiseux polynomial $g:=g(X_0,\dots,X_n)$ is a {\it (tropical) zero} of $f$ (\ref{1}) if for any $(x_0,\dots,x_n)\in \RR^{n+1}$ the value $y=g(x_0,\dots,x_n)$ is a tropical zero of the tropical polynomial $f$.

First, we describe the tropical Newton-Puiseux zeroes of $f$ geometrically and show that there are exactly $n$ of them. In the next section we provide for them the explicit formulas.

For a point $x:=(x_0,\dots,x_n)\in \RR^{n+1}$ its Newton polygon $N_x\subset \RR^2$ is the convex hull of the vertical rays $\{(k,\, c)\, :\, c\ge x_k\},\, 0\le k\le n$. Note that the slopes of the edges of $N_x$ are just the tropical zeroes of $f$.

For a subset $S\subset  \{1,\dots,n-1\}$ consider a convex polyhedron $P_S\subset \RR^{n+1}$ (of dimension $n+1$) consisting of points $x=(x_0,\dots,x_n)$ such that its Newton polygon $N_x$ has the vertices $(0,\, x_0),\, (n,\, x_n),\, \{(s,\, x_s)\, :\, s\in S\}$. Thus,  $\{P_S\, :\, S\subset  \{1,\dots,n-1\}\}$ constitute a partition of $\RR^{n+1}$ into $2^{n-1}$ polyhedra.

Take the (open) polyhedron $P:=P_{\{1,\dots,n-1\}}$ consisting of points $x$ such that the Newton polygon $N_x$ has $n+1$ vertices. Then there are exactly $n$ continuous piece-wise linear functions $g_1,\dots,g_n$ on $P$ being tropical zeroes of $f$ (\ref{1}). Namely, $g_k(x_0,\dots,x_n)=x_{k-1}-x_k,\, 1\le k\le n$.

Observe that each $g_k,\, 1\le k\le n$ has a unique (continuous) continuation on every polyhedron $P_S$. Namely, take the unique pair $0\le i\le k-1,\, k\le j\le n$ such that $i,\, j\in S\cup \{0,\, n\}$, and there are no $s\in S\cup \{0,\, n\}$ satisfying inequalities $i<s<j$.

\begin{lemma}\label{slope}
The unique continuation of $g_k$ on $P_S$ coincides with $\frac{x_i-x_j}{j-i}$.
\end{lemma}

{\bf Proof}. For any point $(x_0,\dots,x_n)$ which belongs to both boundaries of $P$ and of $P_S$ holds $x_s-x_{s+1}=x_{k-1}-x_k,\, i\le s<j$, hence $x_{k-1}-x_k=\frac{x_i-x_j}{j-i}$. $\Box$ \vspace{2mm}

Note that $\frac{x_i-x_j}{j-i}$ is the slope of the edge with the end-points $(i,\, x_i),\, (j,\, x_j)$. Thus, we have shown that there are exactly $n$ tropical Newton-Puiseux polynomials on  $\RR^{n+1}$ being tropical zeroes of $f$ (\ref{1}).

\section{Explicit formulas for tropical zeroes}

\begin{theorem}
A tropical polynomial $f=\min_{0\le k\le n} \{x_k+kY\}$ with parametric coefficients $(x_0,\dots,x_n)$ has exactly $n$ tropical zeroes $g_1,\dots,g_n$ being tropical Newton-Puiseux polynomials in $(x_0,\dots,x_n)$. For each $0\le k\le n$ one can represent $g_k$ as follows. For every $0\le p<k$ consider a tropical Newton-Puiseux polynomial 
$$t_p:=\max_{k\le q\le n} \left\{\frac{x_p-x_q}{q-p}\right\}.$$
\noindent Then $g_k=\min_{0\le p<k} \{t_p\}$. 
\end{theorem}

{\bf Proof}. Fix for the time being a polyhedron $P_S$ and follow the notations from Lemma~\ref{slope}. For any point $x:=(x_0,\dots,x_n)\in P_S$ its Newton polygon $N_x$ has an edge with the end-points $(i,\, x_i),\, (j,\, x_j)$. Therefore, for every $0\le p<k$ the following inequality for the slopes holds:
$$\frac{x_p-x_j}{j-p} \ge \frac{x_i-x_j}{j-i}.$$
\noindent Hence $t_p \ge \frac{x_i-x_j}{j-i}$. 

On the other hand, $t_i=\frac{x_i-x_j}{j-i}$ since for every $k\le q\le n$ the following inequality  for the slopes holds:
$$\frac{x_i-x_q}{q-i} \le \frac{x_i-x_j}{j-i}.$$
\noindent Thus, $g_k$ coincides with $\frac{x_i-x_j}{j-i}$ on $P_S$ which completes the proof due to Lemma~\ref{slope}. $\Box$ 

\begin{remark}
In the formula for $g_k$ in the theorem a tropical Newton-Puiseux polynomial $t_p$ is involved in which a left-end point $(p,\, x_p)$ of the intervals is fixed. In a dual way one can define $r_q:=\min_{0\le p<k} \{\frac{x_p-x_q}{q-p}\}$ by fixing a right end-point $(q,\, x_q)$. Then, similarly to the theorem we get $g_k=\max_{k\le q\le n} \{r_q\}$.  
\end{remark}

\vspace{2mm}
{\bf Acknowledgements}. The author is grateful to the grant RSF 21-11-00283 and
 to MCCME for inspiring atmosphere.

\end{document}